\documentclass[letterpaper, 10 pt, conference]{ieeeconf}  

\IEEEoverridecommandlockouts                              
\usepackage{times}

\usepackage[utf8]{inputenc} 
\usepackage[T1]{fontenc}    
\usepackage{hyperref}       
\usepackage{booktabs}       
\usepackage{amsfonts}       
\usepackage{nicefrac}       
\usepackage{microtype}      
\usepackage{amssymb}
\usepackage{stackengine}
\usepackage{makeidx}
\usepackage{color}
\usepackage{epsfig}
\usepackage{mathtools}
\usepackage{relsize}
\usepackage{multirow}
\usepackage{bm}

\usepackage{microtype}
\usepackage{subfigure}
\usepackage{booktabs} 
\usepackage{stackengine}

\usepackage{algorithm}
\usepackage{algpseudocode}

\usepackage{graphicx}

\algdef{SE}[SUBALG]{Indent}{EndIndent}{}{\algorithmicend\ }%
\algtext*{Indent}
\algtext*{EndIndent}

\newcommand{\icol}[1]{
  \begin{smallmatrix}#1\end{smallmatrix}%
}

\title{\LARGE \bf Online Weight-adaptive Nonlinear Model Predictive Control}

\author{Dimche Kostadinov and Davide Scaramuzza 
\thanks{Dimche Kostadinov and Davide Scaramuzza are with the Robotics and Perception Group ( \url{http://rpg.ifi.uzh.ch} ) at the Dept. of Informatics and Neuroinformatics,~ University of Zurich and ETH Zurich, 8050 Zurich, Switzerland, e-mail: {\tt\small  $\{$dimche, sdavide$\}$@ifi.uzh.ch}}
}

\definecolor{somegray}{rgb}{0.5, 0.5, 0.5}
\newcommand{\darkgrayed}[1]{\textcolor{somegray}{#1}}
\makeatletter
\newcommand*\titleheader[1]{\gdef\@titleheader{#1}}
\AtBeginDocument{%
  \let\st@red@title\@title
  \def\@title{%
    \vskip-3em
    \bgroup\normalfont\large\centering\@titleheader\par\egroup
    \vskip1.5em\st@red@title}
}
\makeatother

\titleheader{\darkgrayed{This paper has been accepted for publication at the \\
IEEE/RSJ International Conference on Intelligent Robots and Systems (IROS), Las Vegas, 2020.}}

\begin{document}


\maketitle

\maketitle
\thispagestyle{empty}
\pagestyle{empty}

\begin{abstract}
Nonlinear Model Predictive Control (NMPC) is a powerful and widely used technique for nonlinear dynamic process control under constraints. In NMPC, the state and control weights of the corresponding state and control costs are commonly selected based on human-expert knowledge, which usually reflects the acceptable stability in practice. Although broadly used, this approach might not be optimal for the execution of a trajectory with the lowest positional error and sufficiently "smooth" changes in the predicted controls. Furthermore, NMPC with an online weight update strategy for fast, agile, and precise unmanned aerial vehicle navigation, has not been studied extensively.  To this end, we propose a novel control problem formulation that allows online updates of the state and control weights. As a solution, we present an algorithm that consists of two alternating stages: (i)  state and command variable prediction and (ii)  weights update. We present a numerical evaluation with a comparison and analysis of different trade-offs for the problem of quadrotor navigation. Our computer simulation results show improvements of up to $\bf 70 \%$ in the accuracy of the executed trajectory compared to the standard solution of NMPC with fixed weights.
\end{abstract}
\section{INTRODUCTION}

In the past, diverse approaches have been used for robust, accurate and fast control \cite{KJAstrom:PID:1995}, \cite{Bellman:1954}, \cite{AlbertoBemporad:2003} and \cite{Houska:2011}. Applied across a wide range of domains, from chemicals to aerospace industries, one of the most powerful and practically useful approaches is NMPC \cite{Allgower:2017}. Its main advantage is that it allows making predictions about the immediate future point under constraints while considering all predicted future points over a given horizon. In recent years, several efficient solutions to  NMPC and improvements have been proposed   \cite{Qin:2000}. One of the most prominent methods for solving NMPC problem is the real-time iteration  (RTI)  scheme \cite{Diehl:2012}. The advantage of RTI is that it allows "on-the-fly" prediction updates: 
as new estimates become available, an iterative solution using only a small number of iteration steps gives the new predictions. 

In the NMPC problem, the state and control weights for the corresponding costs significantly impact the control performance. Usually, fixed weights are carefully selected based on human expert knowledge (platform and trajectory wise tuning) for unmanned aerial vehicle navigation. Under many scenarios, this strategy is preferable. Such weight selection takes into account stability-related properties and exploits the strengths of NMPC. However, whether this approach uses the advantages of the NMPC to the full extent remains an open question.

Additionally, weight tuning might not allow good generalization across a set\footnote{In a sense that a single choice of weights might not allow tracking with the lowest positional error and "smooth" change over the predicted controls.} of diverse trajectories. On the other hand, online and data-adaptive strategies for updating the cost weights during trajectory execution were not studied extensively. In this line, besides the link between NMPC and online learning \cite{Wagener:2019}, other connections to known machine learning paradigms with emphasis on the weigh estimation remain not explored.

To utilize the full control potential of NMPC, we present online weight-adaptive approach. We introduce a novel control problem formulation, where, contrary to using predefined and fixed weights for the state cost, we include the weights of the state cost as a variable in our control problem. This allows us to propose an algorithm that can improve the state and control prediction by optimally updating the weights in an online fashion. Moreover, in our approach, we provide a generalization for a class of weight cost priors and function approximations (including but not limited to neural networks). Also, we give connections of our approach to online learning \cite{Shalev-Shwartz:2014}, reinforcement learning \cite{Sutton:1998}, and metric learning \cite{Kulis:2013}.

\begin{figure}[]
\centering
\begin{center}
\begin{minipage}[b]{.76\linewidth}
\centering
\centerline{\includegraphics[width=\columnwidth]{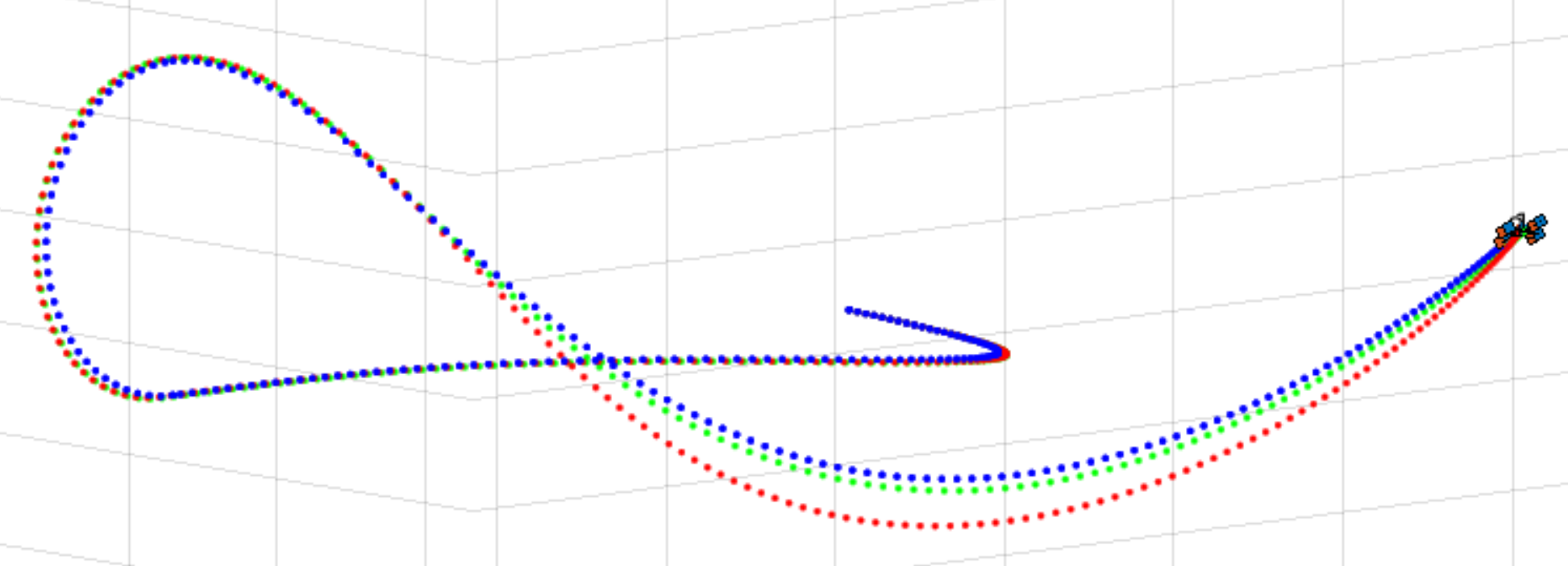}}
\end{minipage}
\begin{minipage}[b]{.73\linewidth}
\centering
\centerline{\includegraphics[width=\columnwidth]{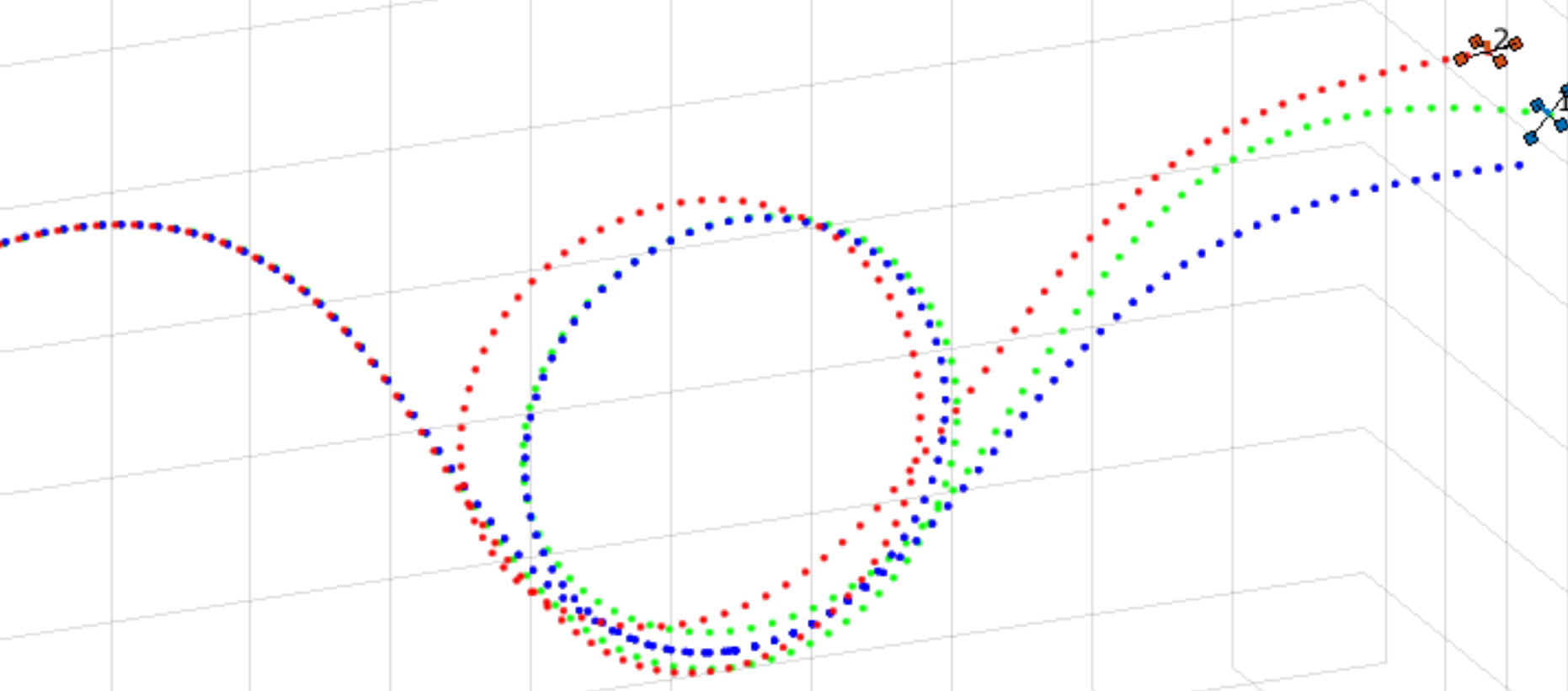}}
\end{minipage}
\caption{{\color{blue} \textbf{A}) The blue line} represents the reference trajectory, {\color{red} \textbf{B}) The red line} represents executed trajectory using the standard solution for NMPC with fixed weights, and {\color{green} \textbf{C}) The green line} represents the executed trajectory using the algorithm that solves the proposed online weight-adaptive NMPC.}
\label{mainMPlot.stage.2}
\end{center}
\vspace{-.1in}
\end{figure}

\subsection{Contributions}

In the following, we summarize our main contributions.

\begin{itemize}
\item We introduce a novel variant of the very well known NMPC control problem, where we address joint prediction of state and control variables and the estimate of the corresponding state and control weights.

\item We propose a  two-stage alternating RTI algorithm. It consists of (i) state and commands variable prediction and (ii) optimal update of the control weights. 

\item We validate our approach by numerical simulation for the task of quadrotor navigation. We demonstrate that our approach reduces the error in the trajectory tracking while predicting controls that have "smooth" change over the executed trajectory. Compared to the solution of the commonly used NMPC, we show improvements of up to $70 \%$ in the accuracy of the executed trajectory. 
\end{itemize}

\subsection{Related Work}

Adaptive Model Predictive Control (MPC) was studied in \cite{An:2009}, \cite{Turki:2017}, \cite{Bujarbaruah:2018}, \cite{Zhao:2017} and \cite{Wagener:2019}. In \cite{An:2009}, the authors proposed an adaptive multi-variable zone controller and gave robustness guarantees for the controlled process. As an adaptive component, they added weighted slack coefficients to the nominal weight coefficients. In \cite{Turki:2017}, the authors were focused on an analytical approach for tuning the control horizon. Their idea consists of computing the value for the optimal control horizon that ensures numerical stability. At the same time, their interest was on a wide set of linear controllable single-input single-output processes. In \cite{Zhao:2017}, the authors developed an adaptive cruise control system for vehicles. The authors utilized a hierarchical control architecture in which a low-level controller compensated for the nonlinear longitudinal vehicle dynamics. Their design enabled tracking of the desired acceleration. They solved a multi-objective control problem by a real-time weight tuning strategy by adjusting each objective's weight for different operating conditions. The authors in \cite{Bujarbaruah:2018}, proposed an adaptive stochastic model predictive control strategy. They considered multi-input multi-output systems that can be expressed by a finite impulse response model. In \cite{Wagener:2019}, a close connection between MPC and online learning was shown. The authors have proposed a new algorithm based on dynamic mirror descent. By doing so, they presented a general family of MPC algorithms that include many existing techniques as special instances. The same authors also provided different  MPC perspectives and suggested principled schemes for the design of new MPC algorithms.

In contrast to the past work, as our main novelty, we consider cost-related weights as an additional variable in the NMPC problem. In this regard, our problem formulation extends the common NMPC problem. We propose a novel algorithm as a solution, while along the way, we also give connections to machine learning paradigms with a focus on metric learning.

\subsection{Paper Organization}

The rest of the paper is organized as follows. In Section \ref{section:Preliminaries}, we present the continuous control problem and its approximate discrete version. In Section \ref{section:ProposedApproach}, we first present the approximate control problem formulation for our NMPC with online weight adaptation. Then, we propose our RTI solution in the form of a two-stage alternating algorithm and discuss the solution. We also give connections to known learning paradigms. While, we devote Section \ref{section:numerical.evaluation} to numerical evaluation, and with Section \ref{section:Conclusions}, we conclude the paper.

\section{ROBOT CONTROL}
\label{section:Preliminaries}
In this section, we first present the continuous problem formulation for robot control and then give its discrete version. 

We assume that the dynamics are described by a set of differential equations $f({\bf x}, {\bf u})$, where ${\bf x} \in \Re^{M_x}$ and ${\bf u} : \Re^{\times M_u}$ denote state and control variables. Furthermore, we assume that an action objective is given, which defines the action cost $\mathcal{L}_a ({\bf x}, {\bf u}): \Re^{M_x \times M_u} \rightarrow \Re^+$. The task of taking action can be expressed by the following  optimization problem:
\begin{equation}
\begin{aligned}
{\{ \hat{\bf x}, \hat{\bf u} \} }=\arg  & \min_{{\bf x}, {\bf u}} \int_{t_0}^{t_h}  \mathcal{L}_a ({\bf x}, {\bf u} )dt,  \\
&\text{ subject to }  f({\bf x} , {\bf u} )=0,  h( {\bf x}, {\bf u} ) \leq 0,
 \label{problem:1}
\end{aligned}
\end{equation}
where $f({\bf x} , {\bf u} )$ and $h({\bf x}, {\bf u} )$ represent equality and inequality constraints that the solution should satisfy when we take action. In order to solve \eqref{problem:1}, first, a cost function for taking action is defined. Then \eqref{problem:1} is discretized, transcribed, and linearized.  In the following text, we go through these steps, which will allow us to present the discretized version of \eqref{problem:1}.

\subsection{The Discrete Control Problem} 

As a common practice, the system dynamics \eqref{problem:1} are  discretized into $N$ system points by using a time step $\delta t$ over fixed time horizon $t_N$. This results in $N$ state vectors ${\bf x}_k$ and $N$ control vectors ${\bf u}_k$ for $k \in \{1,..., N\}$. We assume that we have a reference state ${\bf x}_{r,k}$ and reference controls ${\bf u}_{r,k}$. While, we denote the state errors as $\Delta {\bf x}_{k}={\bf x}_{k}-{\bf x}_{r,k}$ and the control errors as $\Delta{\bf u}_{k}={\bf u}_{k}-{\bf u}_{r,k}$. We define our discrete objective as $\sum_{k=1}^{N}\Delta {\bf x}_{k}^T{\bf Q}_k \Delta {\bf x}_{k}+\sum_{k=1}^{N} \Delta {\bf u}_{k}^T {\bf R}_k \Delta {\bf u}_{k}$, where ${\bf Q}_{k}$ and ${\bf R}_{k}$ denote the state and control weight matrices. 

In \eqref{problem:1}, the equality constraints represents the model for the robot dynamics 
$\frac{ \partial {\bf x}_{k} }{ \partial t}=\dot{\bf x}_{k}=f({\bf x}_k, {\bf u}_k)$. While the inequality constraints $h({\bf x}_k, {\bf u}_k)$ represent the physical limitations of the robot platform. As a common practice, problem  \eqref{problem:1} is transcribed using multiple shooting technique. Moreover, having the discrete dynamics and constraints over the coarse grid $[t_1, ..., t_N]$ for each time interval $[t_{k}, t_{k+1}]$, $\delta t=t_{k+1}-t_{k}$, a boundary value problem is solved, where additional continuity variables are imposed. In addition, an explicit integrator is applied to forward simulate the system dynamics along each interval $[t_{k}, t_{k+1}]$.

Under the above considerations, usually \eqref{problem:1} is sequentially approximated by quadratic problems (QPs). The solution of the QPs are used as an gradient directions $\Delta {\bf x}_k$ and $\Delta {\bf u}_k$ in order to take steps that minimize the original continuous problem. Starting from the available guess (predictions) ${\bf x}_{k}^{pr}$ and ${\bf u}_{k}^{pr}$, the iterations are repeated by taking (not necessarily full) Newton steps in the form $\left[ \icol{{\bf x}^{pr} \\ {\bf u}^{pr}} \right]=\left[ \icol{{\bf x}^{pr} \\ {\bf u}^{pr}} \right]+\alpha \left[ \icol{\Delta {\bf x} \\ \Delta  {\bf u}} \right]$, where $\Delta  {\bf x}=\left[ \icol{ \Delta {\bf x}_1 \\.\\.\\ \Delta  {\bf x}_{N+1}  } \right], \Delta {\bf u}=\left[  \icol{ \Delta  {\bf u}_1 \\.\\.\\ \Delta  {\bf u}_{N}  } \right]$ and  $\alpha$ is the step size which guarantees that the update is in the decent direction. As an example, in the sequential QP approach, given a state estimate ${\bf x}_{1}^m$, a prediction about the state ${\bf x}_{k}^{pr}$ and the control ${\bf u}_{k}^{pr}$, the usual approximation of \eqref{problem:1} with respect to $\Delta {\bf x}_k$ and $\Delta {\bf u}_k$ is the following QP: 
\begin{align}
\! \! \! \! \! \!  \Big[ \! \Delta \hat{\bf x}, \Delta \hat{\bf u} \Big]  \! \! = \! \!   \arg \! \!  \min_{ \Delta {\bf x}, \Delta {\bf u} }      \! \sum_{k=1}^{N} \! & \left( \! \left[\icol{\Delta  {\bf x}_{k} \\ \Delta {\bf u}_{k}}\right]^T \! \! {\bf H}_k \! \! \left[\icol{\Delta {\bf x}_{k} \\ \Delta  {\bf u}_{k}} \right]    \! \! +  \!  \alpha {\bf w}_k^T \left[\icol{\Delta {\bf x}_{k} \\ \Delta  {\bf u}_{k}} \right] \! \right), \nonumber \\
\! \! \! \! \! \! \text{ subject }  \text{to } &  \Delta {\bf x}_{1}={\bf x}_{1}^m-{\bf x}_{1}^{pr}, \label{problem:2.geneeral} \\
& \Delta {\bf x}_{k+1}=  \Delta {\bf r}_k + \left[ \icol{{\bf A}_k, {\bf 0} \\ {\bf 0}, {\bf B}_k} \right] \left[ \icol{\Delta  {\bf x}_{k} \\ \Delta {\bf u}_{k} } \right] ,  \nonumber \\
& \Delta {\bf h}_{k}+\left[ \icol{ {\bf C}_{k}, {\bf 0} \\ {\bf 0}, {\bf D}_{k} } \right] \left[ \icol{ \Delta {\bf x}_{k} \\ \Delta {\bf u}_{k} } \right] \leq {\bf 0}, \nonumber 
\end{align}
where $ \Delta {\bf r}_{k}=f({\bf x}_{k}^{pr}, {\bf u}_{k}^{pr})-{\bf x}_{k+1}^{pr}$, $\Delta {\bf h}_k=h({\bf x}_{k}^{pr}, {\bf u}_{k}^{pr})$. While ${\bf A}_{k}= \frac{\partial f( {\bf x}_k,{\bf u}_k ) }{\partial {\bf x}_k}|_{ \{ {\bf x}_{k}^{pr}, {\bf u}_{k}^{pr} \} }, {\bf B}_{k}= \frac{\partial f( {\bf x}_k,{\bf u}_k ) }{\partial {\bf u}_k}|_{ \{ {\bf x}_{k}^{pr}, {\bf u}_{k}^{pr} \} }$, ${\bf C}_{k}= \frac{\partial h( {\bf x}_k,{\bf u}_k ) }{\partial {\bf x}_k}|_{ \{ {\bf x}_{k}^{pr}, {\bf u}_{k}^{pr} \} }$ and ${\bf D}_{k}= \frac{\partial h( {\bf x}_k,{\bf u}_k ) }{\partial {\bf u}_k}|_{ \{ {\bf x}_{k}^{pr}, {\bf u}_{k}^{pr} \} }$. Where $\frac{\partial f({\bf x}_k, {\bf u}_k)}{\partial {\bf x}_k}|_{\{ {\bf x}_k^{pr}, {\bf u}_k^{pr}\}}$ denoted the partial derivative of $f({\bf x}_k, {\bf u}_k)$ with respect to ${\bf x}_k$, which is evaluated at $\{ {\bf x}_k^{pr}, {\bf u}_k^{pr}\}$,  ${\bf H}_{k}$ is the Hessian of the Lagrangian for \eqref{problem:1} and ${\bf w}_k=\left[ \icol{{\bf Q}_k, {\bf 0} \\ {\bf 0}, {\bf R}_k} \right]{\bf l}_k=\left[ \icol{{\bf Q}_k, {\bf 0} \\ {\bf 0}, {\bf R}_k} \right]\left[ \icol{{\bf x}_k^{pr}-{\bf x}_{r,k} \\ {\bf u}_k^{pr}-{\bf u}_{r,k} } \right]$. One popular approximation to the Hessian ${\bf H}_{k}$ is ${\bf H}_{k}=\left[ \icol{ {\bf Q}_{k}, {\bf 0} \\ {\bf 0}, {\bf R}_k } \right]$ \cite{Sebastien:2016}. We note that in order to simplify the problem description \eqref{problem:2.geneeral}, we omitted the cost related to the last state prediction, but nonetheless, we are taking it into account.

\section{ONLINE WEIGHT-ADAPTIVE NMPC}
\label{section:ProposedApproach}
In this section, we present the problem formulation for our online weight-adaptive  NMPC. Then, we present our two-stage alternating algorithm.  Afterward, we explain and discuss the related problems at each stage. Finally, we give connections to known learning principles. 

\subsection{Min-Max Approximate Control Problem}

We propose to jointly (i) predict the control, and state variables and (ii) estimate the weight matrix. To that end, we present the following problem formulation:
\begin{equation}
\begin{aligned}
\! \! \! \! \! \!  \Big[ \Delta \hat{\bf x},  \Delta \hat{\bf u}, \hat{\bf Q}  \Big] \! = \!  &  \arg   \max_{{\bf Q}}  \bigg\{ \! \! - \lambda g({\bf Q}, \boldsymbol{\theta})+\\ \!  
\! \! \! \! \! \! \min_{ \Delta {\bf x}, \Delta {\bf u} }       \sum_{k=1}^{N} & \left( \left[\icol{\Delta  {\bf x}_{k} \\ \Delta {\bf u}_{k}}\right] + \alpha {\bf l}_k \right)^T  \left[ \icol{ {\bf Q}_{k}, {\bf 0} \\ {\bf 0}, {\bf R}_k } \right]  \left[\icol{\Delta {\bf x}_{k} \\ \Delta  {\bf u}_{k}} \right] \bigg\}, \\
 \! \! \! \! \! \! \text{ subject }  \text{to } & \Delta {\bf x}_{1}={\bf x}_{1}^m-{\bf x}_{1}^{pr},  \\
 & \Delta {\bf x}_{k+1}=  \Delta {\bf r}_k + \left[ \icol{{\bf A}_k, {\bf 0} \\ {\bf 0}, {\bf B}_k} \right] \left[ \icol{\Delta  {\bf x}_{k} \\ \Delta {\bf u}_{k} } \right] ,  \\
& \Delta {\bf h}_{k}+\left[ \icol{ {\bf C}_{k}, {\bf 0} \\ {\bf 0}, {\bf D}_{k} } \right] \left[ \icol{ \Delta {\bf x}_{k} \\ \Delta {\bf u}_{k} } \right] \leq {\bf 0},
 \label{problem:2}
\end{aligned}
\end{equation}
where $\Delta {\bf x}$ and $\Delta {\bf u}$ together with ${\bf Q}= \left[{\bf Q}_{1},..., {\bf Q}_{N} \right]$ are problem variables that we would like to optimally estimate, while $\lambda$ is the Lagrangian variable. In general, $g({\bf Q}, \boldsymbol{\theta})$ could be any weighs related cost function, and $\boldsymbol{\theta}$ is its corresponding parameter. In the simplest form, we consider diagonal ${\bf Q}_k$, and we define $g({\bf Q}, \boldsymbol{\theta})=({\bf Q}{\bf 1})^T({\bf Q}{\bf 1})$, where ${\bf 1}$ is one vector. 

Our proposed formulation is a min-max problem with quadratic and bilinear cost functions. If we fix the weights, the reduced problem over the remaining state and control variables is convex. On the other hand, if we fix the state and control variables then the reduced problem over the weight variables can be  converted again into a convex problem.

\subsection{The Algorithm}

To solve \eqref{problem:2}, we propose an alternating algorithm, which consists of two stages. In stage one, we fix the weights and predict the state and control variables. In stage two, we fix the state and control variables and estimate the weights. In the following, we describe the corresponding problems for the two stages and discuss on their solution.

\subsubsection{\underline{State and Control Prediction}}
Let the weights be fixed, then problem \eqref{problem:2} reduces to the following QP:
\begin{equation}
\begin{aligned}
\! \! \! \! \! \!  \Big[   \Delta \hat{\bf x},  \Delta  \hat{\bf u} \Big]  \! \! =  \! \!   \arg  \! \!  \min_{ \Delta {\bf x} , \Delta {\bf u} } \! &   \sum_{k=1}^{N} \left( \left[\icol{\Delta  {\bf x}_{k} \\ \Delta {\bf u}_{k}}\right]\! \!+\! \alpha {\bf l}_k \right)^T \! \! \left[ \icol{ {\bf Q}_{k}, {\bf 0} \\ {\bf 0}, {\bf R}_k } \right] \! \! \!  \left[\icol{\Delta {\bf x}_{k} \\ \Delta  {\bf u}_{k}} \right], \\
\! \! \! \! \! \!  \text{ subject to } &\Delta {\bf x}_{1}={\bf x}_{1}^m-{\bf x}_{1}^{pr},  \\
 \! \! \! \! \! \! & \Delta {\bf x}_{k+1}= {\bf r}_k + \left[ \icol{{\bf A}_k \\ {\bf B}_k} \right] \left[ \icol{\Delta  {\bf x}_{k} \\ \Delta {\bf s}_{k} } \right],  \\
& \Delta {\bf h}_{k}+\left[ \icol{ {\bf C}_{k} \\ {\bf D}_{k} } \right] \left[ \icol{ \Delta {\bf x}_{k} \\ \Delta {\bf u}_{k} } \right] \leq {\bf 0}.
 \label{problem:2.1}
\end{aligned}
\end{equation}
Since we have quadratic losses, and linear equality and inequality constraints, problem \eqref{problem:2.1} represents a convex quadratic program with linear constraints. Problem \eqref{problem:2.1} is well known and explored \cite{Boyd:2019}, while for a possible solver, we refer to \cite{Boyd:2019}, \cite{Sebastien:2016}. After \eqref{problem:2.1} is solved, as a prediction for the state we use ${\bf x}^{pr}={\bf x}^{pr}+\alpha \Delta \hat{\bf x}$, while as as a prediction for the control we use ${\bf u}^{pr}= {\bf u}^{pr} + \alpha \Delta \hat{\bf u}$.

\subsubsection{\underline{Weights Update}}

This stage enables us to adapt our cost function for taking actions by adjusting and updating the weigh. In turn, this allows us to penalize future errors based on the error between the reference and the currently predicted state. To do so, we fix the control and state variables in \eqref{problem:2} and let ${\bf q}={\bf Q}_1{\bf 1}=...={\bf Q}_N{\bf 1}$. In the simplest form, we define $g({\bf Q})$  as $g({\bf Q}, \boldsymbol{\theta})=({\bf Q}{\bf 1})^T({\bf Q}{\bf 1})$\footnote{It is worthwhile to mention that with $g({\bf Q}, \boldsymbol{\theta})$ we can consider a wide range of parametric function.  Meaning that given data, we could also offline learn and estimate the parameters $\boldsymbol{\theta}$ in our function $g({\bf Q}, \boldsymbol{\theta})$.}. By denoting ${\bf v}_k=\left( \Delta {\bf x}_{k} +\alpha {\bf l}_k \right) \odot \Delta {\bf x}_{k}$,  \eqref{problem:2} reduces to the following quadratic problem:
\begin{equation}
\begin{aligned}
\! \! \! \!  \hat{\bf q}  = 
 & \arg   \min_{{\bf q}}    \lambda{\bf q}^T{\bf q}-\sum_{k=1}^{N_s} {\bf v}_k^T{\bf q}.
 \label{problem:2.2}
\end{aligned}
\end{equation}
where $N_s$, $N_s \leq N$, is the sub-horizon for the weight update. 
The main advantage of \eqref{problem:2.2} is that it has a closed form solution, \textit{i.e.}, $\hat{\bf q}  = \frac{1}{2\lambda+\gamma}\left( \sum_{k=1}^{N_s} {\bf v}_k \right)$. Having the estimated $\hat{\bf q}$, we update ${\bf Q}$ as ${\bf Q}=diag( \hat{\bf q} )$. 

We point out that the vector $\frac{1}{2\lambda+\gamma} \sum_{k=1}^{N_s} {\bf v}_k$ is with non-negative values as long as $\sum_{k=1}^{N_s} \Delta {\bf x}_{k} \geq  \sum_{k=1}^{N_s} \alpha ({\bf x}_{k}^{pr}-{\bf x}_{r,k})$. Therefore, under bounded variations $({\bf x}_{k}^{pr}-{\bf x}_{r,k})$, we can ensure that our weight matrix  ${\bf Q}$ is positive definite. While, under arbitrarily variations a non-negativity\footnote{Note that if we include non-negativity constraint in \eqref{problem:2.2}, the closed from solution for ${\bf q}$ reads as $\hat{\bf q}  = \max \left( \frac{1}{2\lambda+\gamma}\left(  \sum_{k=1}^{N_s} {\bf v}_k\right), {\bf 0} \right)$.} constraint in \eqref{problem:2.2} could be used to ensure that ${\bf Q}$ is positive definite.

\begin{figure}[]
\centering
\begin{center}
\centering
\centerline{\includegraphics[width=\columnwidth]{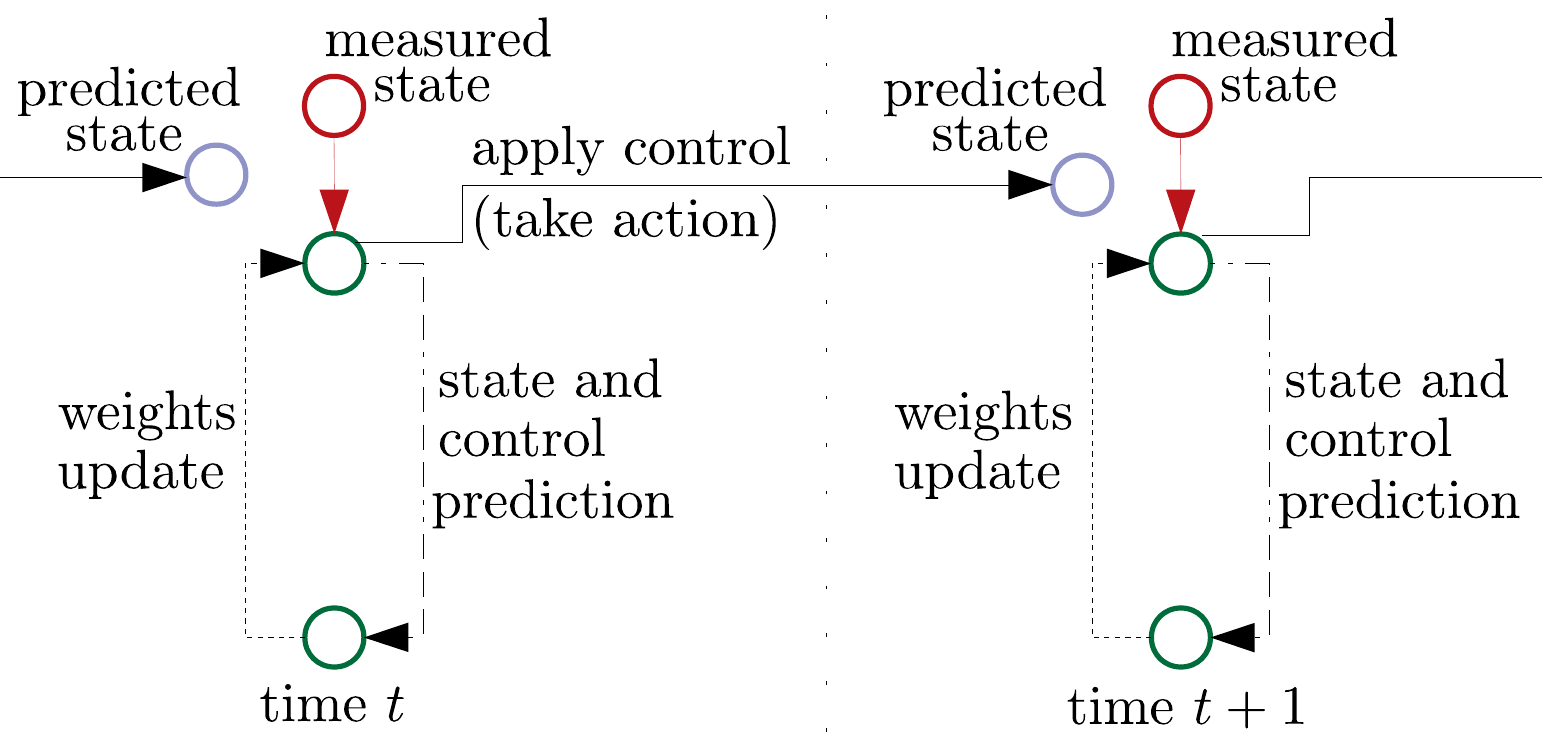}}
\caption{Schematic diagram for the execution stages in our algorithm. In the first stage, we predict the state and control variables. In the second stage, we update the weights. After a number of alternating steps between the two stages, we take action by applying the first control prediction. }
\vspace{-.1in}
\label{intro.one}
\end{center}
\vspace{-.1in}
\end{figure}

\subsection{Connection to Known Learning Paradigms}

Our online weight update approach also represents a special case of metric learning \cite{Kulis:2013}, wherein our case, the loss $\mathcal{L}_a ({\bf x}_k, {\bf u}_k )$ represents the metric, which is included in \eqref{problem:2}. In the following, we give connections to metric learning, online learning, and reinforcement learning.

\subsubsection{Metric Learning} The defined cost function for predicting the state is is the general Mahalanobis distance metric $\left( {\bf x}-{\bf x}_r \right)^T{\bf Q}\left( {\bf x}-{\bf x}_r \right)$. The online update of our distance metric under ${\bf Q}={\bf L}^T{\bf L}$ with ${\bf L}$ is equivalent to learning a linear mapping (given by ${\bf x} \rightarrow {\bf L}{\bf x}$) that transforms the data in the space of ${\bf L}$. After projecting the data onto the new space through the linear map ${\bf L}$, the corresponding distance metric is the usual Euclidean distance.

\subsubsection{Online Learning}
Our algorithm can be viewed as an extension of the online learning approach to model predictive control \cite{Wagener:2019}. Online learning concerns iterative interactions between a learner and an environment over several rounds $N$. At round (or time instance) $t$, the learner picks one out of the set of decisions. The environment then evaluates a loss function based on the learner’s decision, and the learner suffers a cost. The learner’s goal is to minimize the accumulated costs. As shown in \cite{Wagener:2019}, at time t (\textit{i.e.}, round t), an MPC algorithm optimizes a controller (\textit{i.e.}, the decision) over some cost function (\textit{i.e.}, the per-round loss). In this regard, we highlight the following connection to online adaptation and learning. Our alternating algorithm, observes the cost of the initial controller and then improves the controller by updating the cost and the controller, and only then executes a control based on the improved controller.

\begin{algorithm}[t!]
   \caption{\textbf{Online Weight-Adaptive NMPC} 
   }
   \label{alg:example.1}
\begin{algorithmic}
	\Repeat        	
	\State \text{Execute }{\bfseries Stage ${\bf 1}$} 
	\State \text{Execute }{\bfseries Stage ${\bf 2}$} 
    \Until{$convergence$}
   \State {\bfseries Stage ${\bf 1}$} 
   \State \text{ } \text{ } {\bfseries Input } ${\bf x}_1^{m}, {\bf x}_{r}$, ${\bf x}^{pr}$, ${\bf u}_{r}$, ${\bf u}^{pr}$ \text{ and } ${\bf Q}$ 
    \State \text{ } \text{ } \text{ } \text{ } \text{ } $\left[ \icol{\Delta {\bf x} \\ \Delta  {\bf u} } \right] \leftarrow $
    \text{updateDirection}$({\bf x}_1^{m}, {\bf x}_{r}$, ${\bf x}^{pr}, {\bf u}_{r}, {\bf u}^{pr}, {\bf Q})$ 
    
    \State \text{ } \text{ } \text{ } \text{ } \text{ }$\left[ \icol{{\bf x}^{pr} \\ {\bf u}^{pr}} \right]=\left[ \icol{{\bf x}^{pr} \\ {\bf u}^{pr}} \right]+\alpha \left[ \icol{\Delta {\bf x} \\ \Delta  {\bf u} } \right]$
   \State \text{ } \text{ } {\bfseries Output } ${\bf x}^{pr}$ and ${\bf u}^{pr}$
   \State {\bfseries Stage ${\bf 2}$} 
   \State \text{ } \text{ } {\bfseries Input } ${\bf x}^{pr},{\bf u}^{pr},{\bf x}_r$ and ${\bf u}_r$
\State \text{ } \text{ }\text{ } \text{ } \text{ } $ {\bf Q} \leftarrow $
    \text{updateWeight}$({\bf x}_{r}, {\bf u}_{r}$, ${\bf x}^{pr}, {\bf u}^{pr})$ 
   \State \text{ } \text{ } {\bfseries Output } ${\bf Q}$
   
\end{algorithmic}

\end{algorithm}

\subsubsection{Reinforcement Learning}
Regarding the connection of our algorithm to the core principle behind reinforcement learning, we have the following. Upon observing a measurement, in the first stage of our algorithm, we generate state and control prediction. In light of reinforcement learning, we consider this as a sample from some underlining control policy. Afterward, in stage two, we update the weights. Thus our cost metric translates into updating the policy after observing the error between the predicted and reference state.

\section{NUMERICAL EVALUATION}
\label{section:numerical.evaluation}

In this section, we evaluate our approach. Our numerical experiments consider trajectory execution for a quadrotor. Therefore, in the following subsection, we first present the used dynamical model for quadrotor control. Afterward, we describe the experimental setup and discuss the results.

\begin{figure*}[t!]
\centering
\begin{center}
\begin{minipage}[b]{.25\linewidth}
\centering
\centerline{\includegraphics[width=\columnwidth]{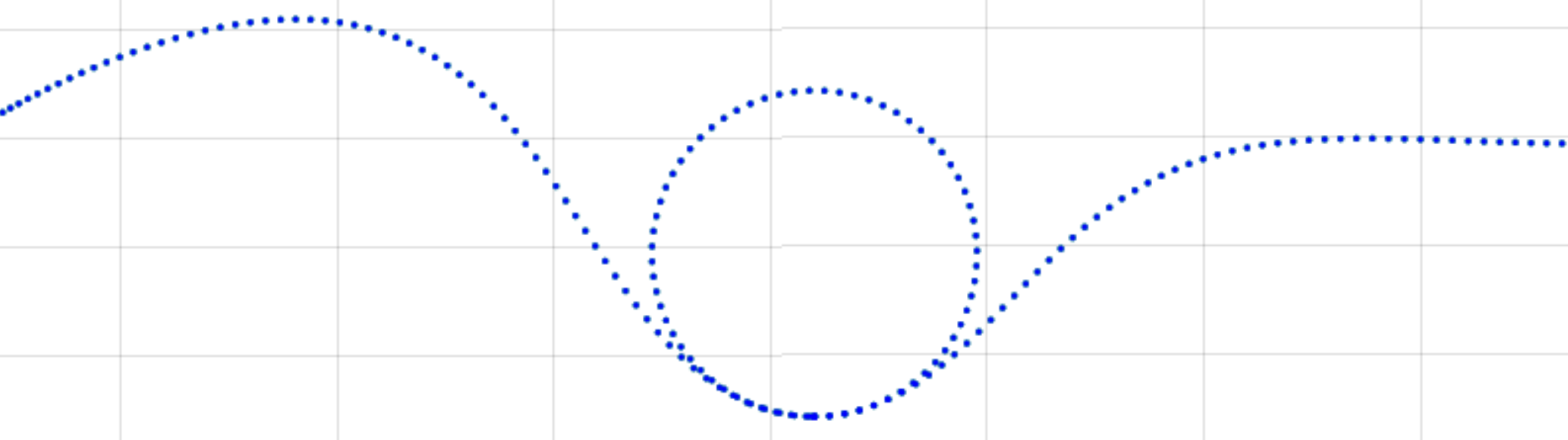}}
{\tiny  \text{Aggressive $T_1$} }
\end{minipage}
\begin{minipage}[b]{.05\linewidth}
\centering
\text{ }
\end{minipage}
\begin{minipage}[b]{.25\linewidth}
\centering
\centerline{\includegraphics[width=\columnwidth]{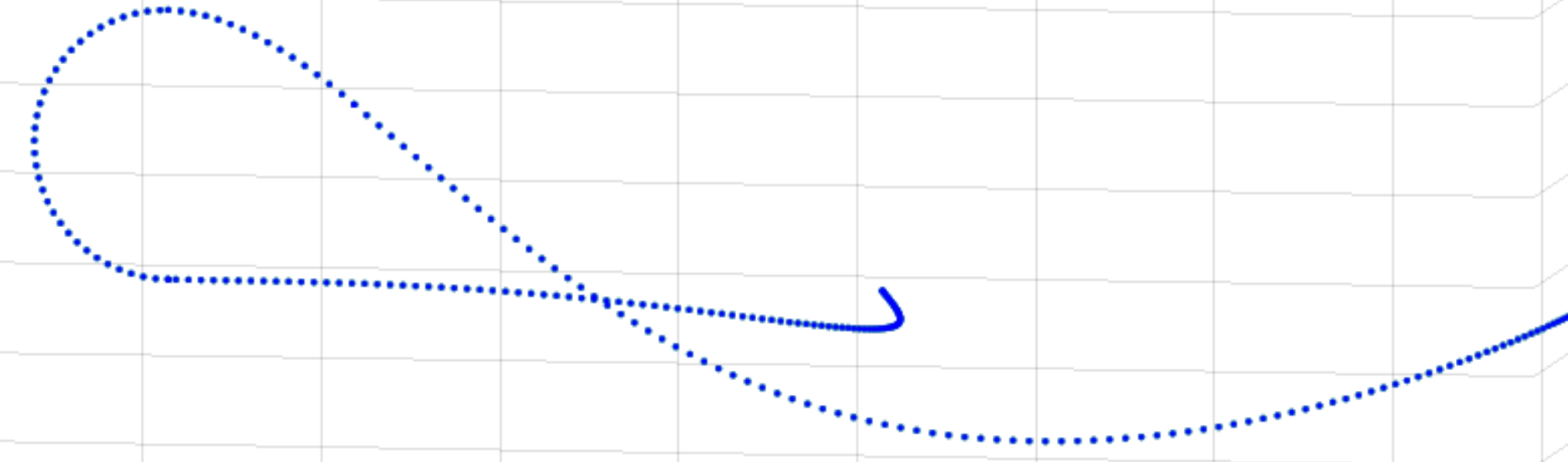}}
{\tiny  \text{Aggressive $T_2$} }
\end{minipage}
\begin{minipage}[b]{.05\linewidth}
\centering
\text{ }
\end{minipage}
\begin{minipage}[b]{.1\linewidth}
\centering
\centerline{\includegraphics[width=\columnwidth]{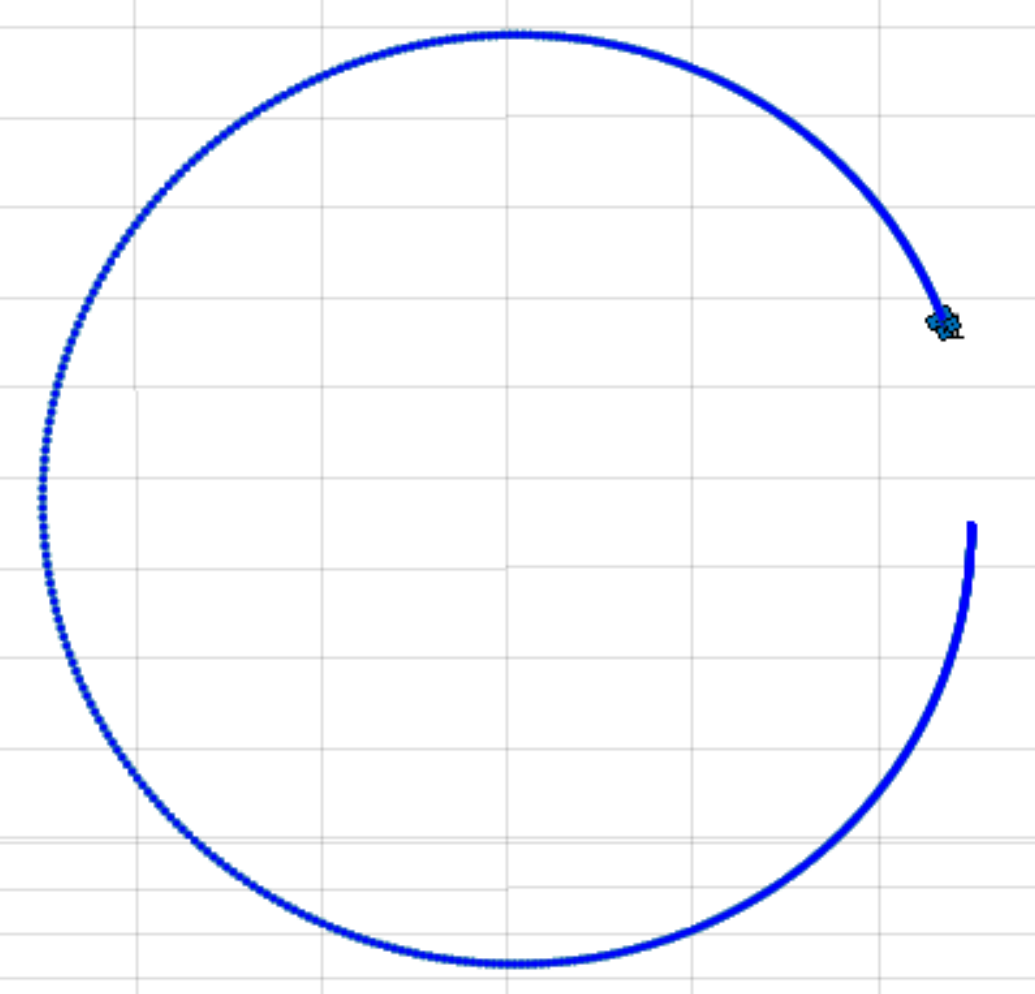}}
{\tiny  \text{Circle $T_3$} }
\end{minipage}
\begin{minipage}[b]{.05\linewidth}
\centering
\text{ }
\end{minipage}
\begin{minipage}[b]{.1\linewidth}
\centering
\centerline{\includegraphics[width=\columnwidth]{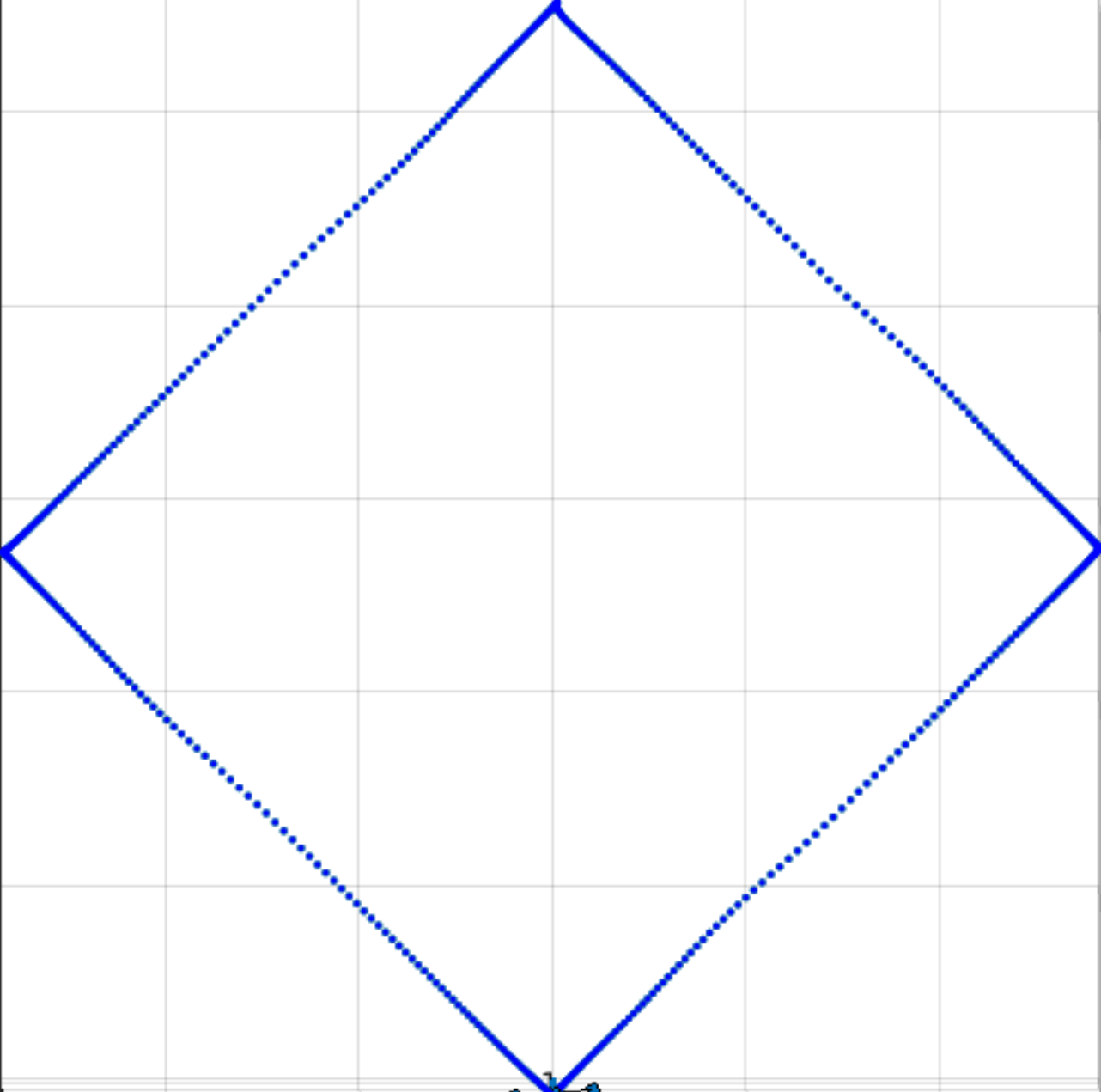}}
{\tiny  \text{Diamond $T_4$} }
\end{minipage}
\caption{Visualization of the used trajectories.} 
\label{trajectorires.visialization}
\end{center}
\vspace{-.1in}
\end{figure*}

\begin{table}[t]
\caption{The total position error $e$ and the cumulative total variation $TV$ over the predicted controls  in the execution of the trajectories under varying $\lambda$.}
\vspace{-.1in}
\label{table.1}
\begin{center}
\begin{tabular}{c|c|c|c|c|c|c|c|}

\cline{3-6}  
\multicolumn{2}{c}{ } & \multicolumn{4}{|c|}{$\lambda$} \\
\cline{3-6}  
\multicolumn{2}{c|}{ } & 0.01   &  0.67  &   1.67   &  3.00 \\
\cline{3-6}  
\multicolumn{6}{c}{ } \\
\cline{3-6}  
\multicolumn{2}{c|}{ } & e[m] | TV   & e[m] | TV  &  e[m] | TV   &  e[m] | TV \\
\cline{3-6}  
\multicolumn{6}{c}{ } \\
\cline{2-6}
\multirow{4}{1em}{\centering $T_1$} 
&\cite{Houska:2011}   &\textbf{0.44} | \textbf{1.57} & \textbf{0.66} | \textbf{1.29} & 1.9 | 1.21 & \textbf{0.87} | \textbf{1.19} \\  
\cline{2-6}
&$N_{2}$  &0.76 | 2.25 & 0.95 | 1.58 & 1.55 | 1.56 & 1.93 | 1.58 \\ 
&$N_{8}$  &0.78 | 2.26 & 1.19 | 1.52 & \textbf{1.36} | \textbf{1.43} & {1.70} | 1.38 \\ 
&$N_{12}$  &0.79 | 1.89 & 1.39 | 1.48 & 1.58 | \textbf{1.43} & 1.82 | 1.41 \\ 
\cline{2-6}
\end{tabular}
\begin{tabular}{cccccccc}
&  &  &  &  & & & 
\hspace{-.7in}
\end{tabular}
\begin{tabular}{c|c|c|c|c|c|c|c|}
\cline{2-6}
\multirow{4}{1em}{\centering $T_2$} 
&\cite{Houska:2011}   &\textbf{0.37} | \textbf{0.72} & 0.64 | \textbf{0.67} & 1.16 | \textbf{0.64} & 1.85 | \textbf{0.63} \\  
\cline{2-6} 
&$N_{2}$  &0.49 | 0.77 & \textbf{0.58} | 0.75 & \textbf{0.68} | 0.69 & \textbf{0.80} | 0.67 \\ 
&$N_{8}$  &0.5 | 0.73 & 0.60 | 0.73 & 0.72 | 0.68 & 0.84 | 0.66 \\ 
&$N_{12}$  &0.52 | \textbf{0.72} & 0.62 | 0.72 & 0.75 | 0.68 & 0.88 | 0.65 \\ 
\cline{2-6}
\end{tabular}
\begin{tabular}{cccccccc}
&  &  &  &  & & & 
\hspace{-.7in}
\end{tabular}
\begin{tabular}{c|c|c|c|c|c|c|c|}
\cline{2-6}
\multirow{4}{1em}{\centering $T_3$} 
&\cite{Houska:2011}   &\textbf{0.34} | \textbf{0.03} & 0.60 | \textbf{0.02} & 2.11 | \textbf{0.02} & 4.49 | \textbf{0.02} \\  
\cline{2-6}
&$N_{2}$  &1.22 | 0.38 & 0.98 | 0.04 & 1.89 | 0.04 & 3.87 | 0.03 \\ 
&$N_{8}$  &0.94 | 0.20 & 0.73 | 0.04 & 1.82 | 0.03 & 3.81 | 0.03 \\ 
&$N_{12}$  &0.68 | 0.86 & \textbf{0.57} | 0.04 & \textbf{1.76} | 0.03 & \textbf{3.72} | 0.03 \\ 
\cline{2-6}
\end{tabular}
\begin{tabular}{cccccccc}
&  &  &  &  & & & 
\hspace{-.7in}
\end{tabular}
\begin{tabular}{c|c|c|c|c|c|c|c|}
\cline{2-6}
\multirow{4}{1em}{\centering $T_4$} 
&\cite{Houska:2011}   &\textbf{0.27} | \textbf{0.08} & 4.63 | \textbf{0.05} & 14.6 | \textbf{0.05} & 28.7 | \textbf{0.04} \\  
\cline{2-6}
&$N_{2}$  &0.38 | 0.21 & 2.56 | 0.12 & \textbf{6.42} | 0.11 & 12.4 | 0.28 \\ 
&$N_{8}$  &0.29 | 0.11 & 2.43 | 0.08 & 6.33 | 0.12 & 12.2 | 0.19 \\ 
&$N_{12}$  &\textbf{0.27} | 0.46 & \textbf{2.03} | 0.46 & 5.87 | 0.14 & \textbf{11.9} | 0.31 \\ 
\cline{2-6}
\end{tabular}
\end{center}
\vspace{-.1in}
\end{table}

\subsection{Used Model}

In the following, we describe the used dynamical model. 

\subsubsection{Quadrotor Dynamics}

Our state ${\bf x}$ and control ${\bf u}$ vectors of the quadrotor are  defined as ${\bf x} = \left[ \icol{ {\bf p}_{WB} \\ {\bf v}_{WB} \\ {\bf q}_{WB} } \right] \text{ $ $} \text{ $ $} \text{ $ $}\text{ and  $ $}\text{ $ $}\text{ $ $} {\bf u} = \left[ \icol{ c \\ \boldsymbol{\omega}_B } \right]$, where ${\bf p}_{WB} = [p_x, p_y, p_z]^T$ and ${\bf q}_{WB} = [q_w, q_x, q_y, q_z]^T$ denote the position and the orientation of the body frame $B$ with respect to the world frame $W$, expressed in world frame, respectively. While ${\bf v}_{WB} = [v_x, v_y, v_z]^T$ denotes the linear velocity of the body, expressed in world frame, and $\boldsymbol{\omega}_B = [\omega_x, \omega_y, \omega_z]^T$ its angular velocity, expressed in the body frame. The vector ${\bf c} = [0, 0, c]^T$ is the mass-normalized thrust vector, where $c = (f_1 + f_2 + f_3 + f_4) /m$, $f_i$ is the thrust produced by the $i$-th motor, and $m$ is the mass of the vehicle. We define the dynamical model for the quadrotor as $\frac{\partial {\bf x}}{\partial t}=f({\bf x}, {\bf u})=\left[ \icol{ \frac{\partial {\bf p}_{WB} }{\partial t}  \\ \frac{\partial {\bf v}_{WB} }{\partial t} \\  \frac{\partial {\bf q}_{WB} }{\partial t}   } \right] = \left[ \icol{ {\bf v}_{WB} \\  _W {\bf g} + {\bf q}_{WB}\odot {\bf c}  \\ \frac{1}{2} \boldsymbol{\Lambda}(\boldsymbol{\omega}_B )  {\bf q}_{WB} } \right]$, where $\frac{\partial {\bf p}_{WB} }{\partial t}$, $\frac{\partial {\bf v}_{WB} }{\partial t}$  and $\frac{\partial {\bf q}_{WB} }{\partial t}$ are the time derivatives of the position, liner velocity and  the quaternion, while  $_W {\bf g}=[0, 0, -g]$ is the gravity vector, with $g=9.81 m/s$. The operator $\odot$ denotes the multiplication between a quaternion and a vector, $\boldsymbol{\Lambda}(\boldsymbol{\omega}_B )  {\bf q}_{WB}$   denotes the time derivative of a quaternion ${\bf q}_{WB}$, while $\boldsymbol{\Lambda}(\boldsymbol{\omega}_B )$ is the the skew-symmetric matrix of the vector $\boldsymbol{\omega}_B$. 

\begin{table}[t!]
\caption{The total position error $e$ and the cumulative total variation $TV$ over the predicted controls in the  execution of the trajectories over varying horizon length $N$.}
\vspace{-.1in}
\label{table.2}
\begin{center}
\begin{tabular}{c|c|c|c|c|c|c|c|}
\cline{3-6}  
\multicolumn{2}{c}{ } & \multicolumn{4}{|c|}{$N$} \\
\cline{3-6}  
\multicolumn{2}{c|}{ } & 8 &   14    &19  &  24    \\
\cline{3-6}  
\multicolumn{6}{c}{ } \\
\cline{3-6}  
\multicolumn{2}{c|}{ } & e[m] | TV   & e[m] | TV  &  e[m] | TV   &  e[m] | TV \\
\cline{3-6}  
\multicolumn{6}{c}{ } \\
\cline{2-6}
\multirow{4}{1em}{\centering $T_1$} 
&\cite{Houska:2011}   &12.0 | \textbf{1.24} & 2.22 | \textbf{1.38} & \textbf{0.98} | \textbf{1.29} & \textbf{0.64} | \textbf{1.34} \\  
\cline{2-6}
&$N_{8}$  &\textbf{1.49} | 1.38 & \textbf{1.38} | 1.49 & 1.13 | 1.56 & 1.22 | 1.55 \\ 
&$N_{14}$  &0 | 0 & 1.61 | 1.58 & 1.45 | 1.53 & 1.22 | 1.55 \\ 
&$N_{18}$  &0 | 0 & 0 | 0 & 1.28 | 1.55 & 1.00 | 1.56 \\ 
\cline{2-6}
\end{tabular}
\begin{tabular}{cccccccc}
&  &  &  &  & & & \\
\end{tabular}
\begin{tabular}{c|c|c|c|c|c|c|c|}
\cline{2-6}
\multirow{4}{1em}{\centering $T_2$} 
&\cite{Houska:2011}   &3.03 | 0.64 & 1.19 | \textbf{0.66} & 0.71 | \textbf{0.67} & \textbf{0.47} | \textbf{0.68} \\  
\cline{2-6} 
&$N_{8}$  &\textbf{0.97} | 0.69 & 0.73 | 0.71 & \textbf{0.62} | 0.72 & 0.53 | 0.74 \\ 
&$N_{14}$  &0 | 0 & \textbf{0.75} | 0.70 & 0.64 | 0.72 & 0.55 | 0.73 \\ 
&$N_{18}$  &0 | 0 & 0 | 0 & 0.65 | 0.72 & 0.56 | 0.73 \\ 
 \cline{2-6}
\end{tabular}
\begin{tabular}{cccccccc}
&  &  &  &  & & & \\
\end{tabular}
\begin{tabular}{c|c|c|c|c|c|c|c|}
\cline{2-6}
\multirow{4}{1em}{\centering $T_3$} 
&\cite{Houska:2011}   &12.75 | \textbf{0.03} & 2.73 | \textbf{0.02} & 0.72 | \textbf{0.02} & \textbf{0.45} | \textbf{0.02} \\  
\cline{2-6}
&$N_{8}$  &\textbf{4.05} | 0.05 & 1.79 | 0.04 & 0.86 | 0.04 & 0.52 | 0.04 \\ 
&$N_{14}$  &0 | 0 & \textbf{1.17} | 0.07 & 0.66 | 0.04 & 0.50 | 0.04 \\ 
&$N_{18}$  &0 | 0 & 0 | 0 & \textbf{0.61} | 0.92 & 0.51 | 0.10 \\ 
\cline{2-6}
\end{tabular}
\begin{tabular}{cccccccc}
&  &  &  &  & & & \\
\end{tabular}
\begin{tabular}{c|c|c|c|c|c|c|c|}
\cline{2-6}
\multirow{4}{1em}{\centering $T_4$} 
&\cite{Houska:2011}   &59.62 | 0.06 & 15.7 | \textbf{0.05} & 5.73 | \textbf{0.05} & 2.07 | \textbf{0.05} \\  
\cline{2-6}
&$N_{8}$  &\textbf{24.69} | 0.07 & 6.90 | 0.08 & 2.94 | 0.08 & 1.44 | 0.08 \\ 
&$N_{14}$  &0 | 0 & \textbf{5.48} | 0.08 & 2.27 | 0.56 & 1.19 | 0.42 \\ 
&$N_{18}$  &0 | 0 & 0 | 0 & \textbf{1.95} | 0.74 & \textbf{1.15} | 0.73 \\ 
\cline{2-6}
\end{tabular}
\end{center}
\vspace{-.1in}
\end{table}

\subsubsection{Quadrotor Physical Constraints} By the inequality constraint $ h({\bf x}, {\bf u})$, we model the physical limitations of the drone platform in order to attain feasible solutions. In our case it is the minimum and maximum thrust $c_{min}$ and $c_{max}$, as well as minimum and maximum angular velocities $\boldsymbol{\omega}_{min}$ and $\boldsymbol{\omega}_{max}$, respectively, which we compactly express as ${\bf u}_{min}=[\icol{c_{min}\\ \boldsymbol{\omega}_{min}}]$ and ${\bf u}_{max}=[\icol{c_{max}\\ \boldsymbol{\omega}_{max}}]$.

\subsection{Setup, Error Measures and Comparative Analysis}

We generated  four different trajectories, which were computed as proposed in \cite{Mellinger:2011}. As shown in Figure \ref{trajectorires.visialization}, the trajectories have different geometries. The first and the second trajectory are aggressive and are denoted as $T_1$ and $T_2$. The third trajectory $T_3$ is circle and the fourth trajectory $T_4$ is diamond. 

We validated our approach under different setups. We experimented with different strength for the state costs in the control  problem. As well as we experimented with  different lengths of the fixed horizon and different lengths of the sub horizon, which were used to update the weight in the cost. In addition, we validated our approach under additive white Gaussian noise perturbation in the available state estimate. In summary, we present simulation results under:
\begin{itemize}
\item[(i)] Different strength $\lambda$ of the state cost,
\item[(ii)] Different length of the prediction horizon $N$ and 
\item[(iii)] Noise perturbation with different noise levels $\sigma$ in the available (measured) state.
\end{itemize}

\begin{table}[t]
\caption{The total position error $e$ and the cumulative total variation $TV$ over the predicted controls in the  execution of the trajectories under noise perturbation with noise level $\sigma$.}
\vspace{-.1in}
\label{table.3}
\begin{center}
\begin{tabular}{c|c|c|c|c|c|c|c|}
\cline{3-6}  
\multicolumn{2}{c}{ } & \multicolumn{4}{|c|}{$\sigma$} \\
\cline{3-6}  
\multicolumn{2}{c|}{ } & 0.5 &   2.0    & 3.5  &  5.0    \\
\cline{3-6}  
\multicolumn{6}{c}{ } \\
\cline{3-6}  
\multicolumn{2}{c|}{ } & e$_r[m]$   & e$_r[m]$  &  e$_r[m]$   &  e$_r[m]$ \\
\cline{3-6}  
\multicolumn{5}{c}{ }  \\
\cline{2-6}
\multirow{2}{1em}{\centering $T_1$} 
&\cite{Houska:2011}   &\textbf{1.32} & \textbf{4.48} & 8.04 & 12.32 \\  
\cline{2-6}
&$N_{8}$  &1.65 & 4.55 & \textbf{7.58} & \textbf{11.58} \\ 
\cline{2-6}
\multicolumn{6}{c}{ } \\
\cline{2-6}
\multirow{2}{1em}{\centering $T_2$} 
&\cite{Houska:2011}   &1.25 & 4.79 & 8.23 & 12.48 \\  
\cline{2-6}
&$N_{8}$  &\textbf{1.13} & \textbf{4.31} & \textbf{7.46} & \textbf{11.21} \\ 
\cline{2-6}
\multicolumn{6}{c}{ } \\
\cline{2-6}
\multirow{2}{1em}{\centering $T_3$} 
&\cite{Houska:2011}   &1.77 & 5.38 & 9.04 & 12.87 \\  
\cline{2-6} 
&$N_{8}$  &\textbf{1.44} & \textbf{4.87} & \textbf{8.95} & \textbf{12.57} \\ 
\cline{2-6}
\multicolumn{6}{c}{ } \\
\cline{2-6}
\multirow{2}{1em}{\centering $T_4$} 
&\cite{Houska:2011}   &3.95 &    5.07  &  \textbf{8.72} & 12.74 \\  
\cline{2-6}
&$N_{8}$  &\textbf{1.3} & \textbf{4.8} & 8.85 & \textbf{12.18} \\ 
\cline{2-6}
\end{tabular}
\end{center}
\vspace{-.1in}
\end{table}

Over all trajectory points, we measure and report the total accumulation of error as  $e=\sum_{i=1}^L d_i $, where $d_i=\Vert \left[\icol{ p_x \\ p_y \\ p_z  }\right]-\left[\icol{ p_{r,x} \\ p_{r,y} \\ p_{r,z}  } \right] \Vert_2$ represents the error between the simulated position $\left[  \icol{ p_{x} \\ p_{y} \\ p_{z}  } \right]$ after applying the predicted control and the reference position $\left[ \icol{ p_{r,x} \\ p_{r,y} \\ p_{r,z} }\right]$. In addition, we also measure the total variation over commands, \textit{i.e.}, $\text{TV} = \frac{1}{L}\sum_{i=1}^{L-1}\left( \vert T_i-T_{i+1} \vert + \sum_{j=1}^3\vert \omega_{j,i}-\omega_{j,i+1} \vert  \right)$ and consider it as an indicator for "smooth" changes in the predicted commands during trajectory execution.  

In the first two series of experiments for different strength $\lambda$ of the state cost, and for different lengths of the prediction horizon $N$, we also experimented with different sub-horizon lengths $N_s$.  Regarding the noise corruption protocol for the third set of experiments, we implement it as follows. We randomly selected the index $\tau \in \{1,..., L \}$ for one point from the trajectory. Then, during execution, at the corresponding index $\tau$, we corrupted the available state with  Additive White Gaussian Noise (AWGN) $\boldsymbol{\nu}$ as follows ${\bf p}_{\tau}^{'}={\bf p}_{\tau}+\sigma \boldsymbol{\nu}$, while we ensured that $\Vert {\bf p}_{\tau} \Vert_2 = \Vert \boldsymbol{\nu} \Vert_2 $. In the same experiment, we compute average error as  $e_r=\frac{1}{K}\sum_{k=1}^{K}e_k$  for $K=15$ runs of this procedure. 
We compare our algorithm with the standard solution to NMPC with fixed weights, which we implemented using \cite {Houska:2011}. 

In the standard NMPC, the hole weight matrix ${\bf Q}$ has to be tuned. In contrast, in our approach the only tuning parameter is $\lambda$. During simulation, we also found out that $\hat{\bf q}  =\exp \Big(  \frac{1}{2\lambda+\gamma}\left( \sum_{k=1}^{N_s} {\bf v}_k \right) \Big)$, has better performance then $\hat{\bf q}  = \frac{1}{2\lambda+\gamma}\left( \sum_{k=1}^{N_s} {\bf v}_k \right)$.  Therefore, in all of the experiments, we update the weights as  $\hat{\bf q}  =\exp \Big(  \frac{1}{2\lambda+\gamma}\left( \sum_{k=1}^{N_s} {\bf v}_k \right) \Big)$.

\subsection{Results Discussion}

In Tables \ref{table.1}, \ref{table.2} and  \ref{table.3}, we present the results of our computer simulation. We show the resulting  trajectory tracking errors and total variation errors of our algorithm and the comparing method (the standard NMPC, with fixed weights). 

As we can see in Table \ref{table.1}, the total error in position is reduced as a result of a small increase in the total variation of the command prediction when compared to the common solution \cite{Houska:2011} of the NMPC. Moreover, as shown in  Table \ref{table.1}, for very small values of the $\lambda$ parameter, the accuracy for the trajectory tracking of our algorithm is lower than the accuracy of the comparing algorithm \cite{Houska:2011}. However, we note that at such small values for $\lambda$, the changes in the predicted controls during the trajectory execution are not "smooth". In practice, very small $\lambda$ might lead to potentially non-stable behavior. On the other hand for $\lambda$ values above $.5$, \textit{i.e.}, $\lambda >.5$, the changes in the predicted controls are more "smooth".  At the same $\lambda$ values, we report improved performance. Our approach achieves lower total error compared to the common solution of the standard NMPC. When the sub-horizon $N_s$\footnote{Our results are computed for sub-horizon length $N_s$ smaller then the horizon length $N$, $N_s \leq N$.} has larger length, the results error is lower, but the TV is also high. 

Table \ref{table.2} shows that the total error in position is consistently lower compared to \cite{Houska:2011} over different horizon lengths. It is interesting to highlight that even for low horizon length like $N=8$, the algorithm achieves high tracking accuracy. While, both of the comparing algorithms have the lowest errors for horizons between $16$ and $24$. 

In Table \ref{table.3}, we can see that for a noise level in the range of $.5$ to $5$, the average positional error $e_r$ of the our algorithm is lower compared to the same error $e_r$ for the standard NMPC \cite{Houska:2011}.


As a summary, the simulation results demonstrated that by our approach, which is without manual weight tuning, we could archive accurate and stable quadrotor navigation. The execution of smooth as well as fast and rapidly changing reference trajectories benefits from online weigh adaptation. The results also show that we can have relatively good tracking performance even under a small horizon length. It is essential to point out that not all online weight update configurations are useful. In other words, not all algorithm setups (for different $\lambda$ and $N_s$) provide improved accuracy with "smooth" changes in the predicted controls over the executed trajectories in the simulation.

\section{CONCLUSIONS}
\label{section:Conclusions}

In this paper, we presented a  novel control problem formulation for NMPC with an online update of the cost weights.  As a  solution,  we proposed a  two-stage alternating algorithm. It consists of: (i) state and commands variable prediction and  (ii)  optimal weights update. Our evaluation by computer simulation demonstrated not only high accuracy for trajectory tracking but also robustness to noise perturbation. Comparing the solution of our approach to the solution of the common NMPC with fixed weights, we demonstrated lower tracking error for the used reference trajectories. Our next steps are to test the performance on a real drone platform.

\bibliography{root}
\bibliographystyle{plain}

\end{document}